\documentclass{article}
\usepackage[hmargin={2truecm,2truecm},vmargin={2truecm,2truecm}]{geometry}

\usepackage{amsmath,amsthm,amsopn,amstext,amscd,amsfonts,amssymb,lscape}
\usepackage{url,color,colortbl}
\usepackage{graphics}
\usepackage{graphicx}
\usepackage{multirow}
\usepackage[normalem]{ulem}
\usepackage{rotating}
\usepackage[ruled,vlined,linesnumbered,noresetcount]{algorithm2e}
\usepackage{float}

\definecolor{Gray}{gray}{0.9}

\newtheorem{propo}{\bf Proposition}[section]
\newtheorem{coro}{\bf Corollary}[section]
\newtheorem{prope}{\bf Property}[section]
\newtheorem{observation}{\bf Observation}[section]

\newcommand{\dem}{\par \noindent{\bf Proof: }}

\newcommand{\MIN}{b-minC}
\newcommand{\MAX}{b-maxC}

\definecolor{arsenic}{rgb}{0.23, 0.27, 0.29}
\definecolor{mygray}{gray}{0.6}
\definecolor{armygreen}{rgb}{0.19, 0.53, 0.43}
\definecolor{atomictangerine}{rgb}{1.0, 0.6, 0.4}

\title{{\sc Budget-constrained cut problems}}

\author{Justo Puerto\footnotemark[1], Jos\'e L. Sainz-Pardo\footnotemark[2]~ \footnotemark[3]~ \footnotemark[4]\\}

\date{}

\begin{document}

\maketitle

\renewcommand{\thefootnote}{\fnsymbol{footnote}}

\footnotetext[1]{IMUS and Departamento de Estadística e Investigación Operativa, Universidad de Sevilla, Spain. E-mail address:puerto@us.es}
\footnotetext[2]{Centro de Investigación Operativa, Universidad Miguel Hernandez de Elche, Spain. E-mail address: jose.sainz-pardo@umh.es}
\footnotetext[3]{Corresponding author, e-mail: jose.sainz-pardo@umh.es}
\footnotetext[4]{This work was written during a sabbatical leave of the second author at IMUS funded by project CIBEST/2021/35 of the Generalitat Valenciana}

\begin{abstract}
The minimum and maximum cuts of an undirected edge-weighted graph are classic problems in graph theory. While the Min-Cut Problem can be solved in $P$, the Max-Cut Problem is NP-Complete. Exact and heuristic methods have been developed for solving them. For both problems, we introduce a natural extension in which cutting an edge induces a cost. Our goal is to find a cut that minimizes the sum of the cut weights but, at the same time, restricts its total cut cost to a given budget. We prove that both restricted problems are NP-Complete and we also study some of its properties. Finally, we develop exact algorithms to solve both as well as a non-exact algorithm for the min-cut case based on a Lagreangean relaxation that generally provides optimal solutions. Their performance is reported by an extensive computational experience. 
\end{abstract}
\medskip

\section{Introduction}

Graph disconnection is a classic topic in graph theory and has been widely studied because it has many applications such as chip and circuit design, network reliability, supply chain planning, cloud and streaming architecture design, and cluster analysis (see e.g. \cite{justo3,Becker_Lari_al_98,benati,justo1,justo2}). Finding the minimum and the maximum cuts of undirected edge-weighted graphs, that is, to find the partition of the graph in two components $S$ and $T$ in such a way that the sum of the edges connecting both components is minimum or maximum, are fundamental problems in graph disconnection. \\

Regarding minimum cuts, there are exact algorithms, such as Stoer-Wagner, that solve the problem in polynomial time. Since the minimum cut of a network is the minimum cut among all the cuts between two nodes $s$ and $t$ such that $s \in S$ and $t \in T$, most of the algorithms compute the global minimum by calculating combinations of $s-t$ cuts, as for instance Stoer-Wagner does. Besides, it is well-known the duality between the minimum $s-t$ cut and the maximum flow between the vertices $s$ and $t$ also being the maximum global flow of a graph the maximum $s-t$ flow among any pair of vertices. Furthermore, compact linear models have also been proposed to determine minimum cuts of a graph.\\

A reduction of the Min-Cut Problem is the Min $s-t$ Cut Problem. In this case one is searching  for a minimum cut in such a way that node $s$ belongs to a cluster and node $t$ to the other one. In fact, most of the algorithms and methods to solve the general Minimum Cut problem resort to solve $|V|-1$ minimum $s-t$ cut problems with different combinations of $s$ and $t$. Then, the minimum cut problem is solved by the minimum of the $s-t$ already found. These approaches encompasses the Stoer-Wagner algorithm as well as linear programming models like the one in \cite{linear} whose constraints can be interpreted as an intersection of $ |V| - 1$ minimum s-t cut polyhedra. Moreover, for solving the general problem, an approach taking advantage of global connectivity has been proposed in  \cite{minCutModel} which results in a linear formulation without integer variables and with complexity $O(|V| |E|)$.

Maximum cut is one of the simplest graph partitioning problems to conceptualize, and yet it is one of the most difficult combinatorial optimisation problems to solve. Contrary to the Min-Cut Problem, which is solvable in polynomial time, the maximum cut problem is NP-complete \cite{complejidad} and is one of Karp's $21$ NP-complete problems, see \cite{reducibilidad}.\\

The inclusion of some additional conditions on the weight or cardinality of the cuts slightly modify the problems but in most cases the addition of constraints is required. These problems are usually called budget-constrained problems. When this is the case, problems easily solvable in polynomial time become NP-Complete: the budget constraint transforms the problem into a version of a knapsack problem. Under this perspective, the constrained-shortest path problem that converts a classic problem like the shortest path, solvable by Dijtstra's Algorithm (\cite{dijkstra}) in $O((E + V) log V)$ time, into a NP-Complete problem, has been well studied. Several works were developed for solving it. Already in 1966, it was proposed an algorithm based on Dynamic Programming \cite{shortest1}. \cite{shortest2} is also based on Dynamic Programming but other approaches are based on Lagrangian relaxations (\cite{shortest3}, \cite{shortest4}, \cite{shortest5}). Although Lagrangian relaxation does not certify optimality of the solutions obtained, the optimal solution can be reached by strengthening with k-shortest enumeration (\cite{kShortest1}, \cite{kShortest2}) . Nevertheless, all these works were surpassed in 2012 by the introduction of the Pulse Algorithm \cite{pulse}, a Dynamic Programming algorithm that despite its simplicity is capable of handling these problems in large-scale networks. Later, the Pulse Algorithm was refined in \cite{pulseBi} by a bidirectional search strategy leveraged on parallelism. \\

Other budget-constrained problems studied, though not to the same extent as the shortest path, are the prize-collecting travelling salesman with a budget-constraint, the budgeted minimum spanning tree problems or the budget constraint flow problem (both in its minimum and maximum versions \cite{flujo3}, \cite{flujo1}, \cite{flujo2}) for which different approaches have been developed (\cite{otros}). However, as far as we know, network cut problems have not been analyzed in their constrained versions despite the usefulness of their applications, although some of them are outlined in \cite{ataque} and \cite{budgetedMax}. In this paper we study the minimum and maximum budget-constrained cut problems. We prove that both are NP-Complete and study some of their properties. We also develop a branch and bound algorithm to solve them that gets results in short computing time for the the Minimum Cut Problem, in the same way that the Pulse algorithm does for the shortest path despite being very different both in concept and purpose. However, in the  case of the Max-Cut Problem the time employed is just reasonable. We also provide two approaches based on Lagrangian relaxation that are capable of providing good approximate solutions of the Min-Cut problem with little computational effort although the optimality of the obtained solutions is not assured. Finally, the performance of the proposed algorithms is validated by a computational experience making a comparative analysis versus the compact integer programming models.\\

\section{The problem}

We consider an undirected graph $G=(V, E)$, with $|V|=n$ and $|E|=m$, such that both the set of vertices $V$ and the set of edges $E$ are associated with non-negative weights $w_v$ and $w_e$, respectively. \\


\subsection{{The Budget Constrained Minimum Cut Problem}}

The Budget-Constrained Min-Cut Problem can be expressed by the integer linear programming formulation (\MIN). 

\begin{alignat}{2}
\mbox{(\MIN)} & \min\ \displaystyle \sum_{(i,j) \in E} w_{ij} x_{ij} \label{fobj2}\\
\mbox{\bf s.t.}\ & &  \nonumber \\
&\sum_{2 \leq i \leq n} z_i=1 & \label{corte1} \\
& x_{ik} + x_{jk} \geq x_{ij} + 2 z_k \hspace{0.5cm} & i, j, k \in V: i<j<k \label{corte2} \\
& x_{ik} \geq z_k & i \in V, k \in V: i<k \label{corte3} \\
& \sum_{(i,j) \in E} c_{ij} x_{ij} \leq T \label{budget} \\
& x_{ij} \in \{0,1\} & i, j \in V \label{binX} \\
& z_i \in \{0,1\}  & i \in V \setminus{1} \label{binZ} \\
\end{alignat}

Objective function and constraints \eqref{corte1}-\eqref{corte3} are borrowed from \cite{minCutModel} which as far as we know is the best compact formulation with $\mathcal{O}(|V|^2)$ constraints and  $\mathcal{O}(|V|^3)$ variables. In the a straightforward case without  budget constraints, the integer character of the variables is not required, since this formulation has integrality property. They use two types of variables $x_{ij}$ and $z_i$. $x_{ij}$ takes value $1$ if edge $(i,j)$ is cut and $z_i$ takes value $1$ if node $i$ is clustered in a group different to the group of node $1$. So, $z_i$ is defined for $i \in \{2, \dots, n\}$. Then, the objective function \eqref{fobj} minimizes  the total weight of the cut. Constraint \eqref{corte1} enforces that there must exist a cluster different from that containing node $1$. In the following constraints, \eqref{corte2}, the node of highest index is considered as representative of this cluster. Therefore, \eqref{corte2} forces variables $x_{ik}$ and $x_{jk}$ to assume value 1 if nodes $i$ and $j$ belong to the cluster of node $1$ but node $k$ does not. Recall that, without loss of generality, we are considering a complete graph. Finally, \eqref{budget} represents the budget constraint. The reader may note in passing that in the original model, $x_{ij}$ and $z_i$ can be relaxed to be continuous. However, when we introduce the budget constraint, \eqref{corte1}-\eqref{corte3} do not represent the convex hull of the solution space, so it is compulsory to require the integrality of these two sets of variables $x_{ij}$ and $z_i$.\\

\subsection{The Budget-Constrained Maximum Cut Problem}

(\MAX) model formulates the Budget-Constrained Maximum Cut Problem as an integer linear programming model.

\begin{alignat}{2} 
\mbox{(\MAX)} & \max\ \displaystyle \sum_{(i,j) \in E} w_{ij} x_{ij} \label{fobj}\\
\mbox{\bf s.t.}\ & &  \nonumber \\
& x_{ik} + x_{jk} \geq x_{ij} \hspace{0.5cm} & i, j, k \in V: i<j<k \label{corteMax1} \\
& x_{ij} + x_{ik} + x_{jk} \leq 2 \hspace{0.5cm} & i, j, k \in V: i<j<k \label{corteMax2} \\
& \sum_{(i,j) \in E} c_{ij} x_{ij} \leq T \label{budgetMax} \\
& x_{ij} \in \{0,1\} & i, j \in V \label{varMax} \\
\end{alignat}

The binary variables $x_{ij}$ have the same meaning as in \MIN\ and its objective function maximizes the total weight of the cut whilst constraints \eqref{corteMax1} and \eqref{corteMax2} have been borrowed again from \cite{max-cut}: constraint \eqref{corteMax1} implies that if edge $(i,j)$ is in the cut, then  edge $(j,k)$ or edge $(i,k)$ must be cut, and constraint \eqref{corteMax2} ensures the triangle inequality. Finally, \eqref{budgetMax} adds the budget constraint. 

\section{Some properties}
\begin{propo} \label{NPCompleto}
The Budget-Constrained Min s-t Cut Problem is NP-Complete.
\dem To prove it, firstly, we are going to transform the Binary Knapsack Problem (BKP) into a budget-constrained Min s-t Cut Problem in such a way that the optimal solutions of both problems coincide. Then, given a set of $n$ items, each with profit $w_j$ and a weight $c_j$ for $j=1,...,n$ we seek to pack them into a knapsack maximizing the profit of the packed items in such a way that the weight limit $T$ is not exceeded. To prove that Budget-Constrained Min-Cut Problem is NP-Complete, we set up a network with $n+2$ nodes: $V=\{s, t, v_1, ..., v_n\}$. On the one hand, each one of the nodes $V \setminus \{s,t\}$ is linked to node $s$ having $0$ weight and $M$ cost. On the other hand, each one of the nodes $v_i \in V \setminus \{s,t\}$ is linked to node $t$ having $w_i$ weight and $M-c_i$ cost. Figure \ref{NPC} illustrates the previously described network.

\begin{figure}[H]
	\begin{center}
		\includegraphics[width=0.3\textwidth]{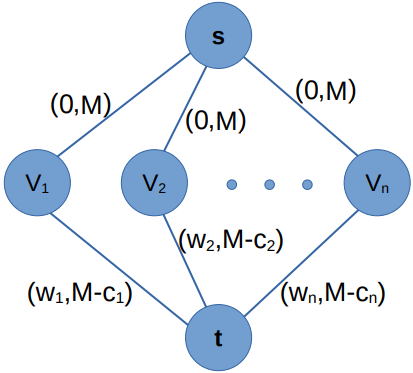}
	\end{center}
	\caption{Graph for the Restricted Min-Cut s-t Problem}\label{NPC}
\end{figure}

BKP can be solved by finding a Budget-Constrained Min s-t Cut Problem on the described graph. Since BKP is NP-Complete, then Budget-Constrained Min s-t Cut Problem is also NP-complete. $\hfill \square$ 

\end{propo}


\begin{propo} \label{minNPCompleto}
The Budget-Constrained Min-Cut Problem is NP-Complete.
\dem As it is well known, the s-t min-cut is a reduction of the global min-cut but if the cuts are restricted it is the same. To show this, we use the proof of Karger and Stein in \cite{karger} for directed graphs. Thus, we add for each vertex v, directed edges of infinite weight and cost $0$ from $t$ to $v$ and from $v$ to $s$. The global minimum cut in this modified graph must now have $s \in S$ and $t \in T$ since otherwise some of the edges of infinite weight will appear in the cut. Since Restricted Min s-t Cut Problem is NP-Complete then budget constrained Min-Cut Problem is NP-Complete. $\hfill \square$

\end{propo}

\begin{propo} \label{maxNPCompleto}
The Budget Constrained Max-Cut Problem is NP-Complete.
\dem Since the Budget Constrained Max-Cut Problem comes from two NP-Complete Problems: the Max-Cut Problem and the Binary Knapsack Problem then the Budget Constrained Max-Cut Problem is NP-Complete. The reader may note that a proof similar to the one in Proposition \ref{minNPCompleto} also holds. $\hfill \square$
\end{propo}

Properties \ref{arbol} and \ref{bloque} are straightforward and the proof is left for the reader.

\begin{prope} \label{arbol} The Budget Constrained Min-Cut (Max-Cut) of a tree is the edge of minimum (maximum) weight such that its cost is lower than the budget.
\end{prope}


In a graph, a component is a maximal connected subgraph. On the one hand, a vertex is called cut-vertex if its removal augments the number of components of the graph. On the other hand, a block of a graph is a maximal connected subgraph without cut-vertices. A separable graph is the one that has cut-vertices. The subgraphs that are connected by cut-vertices are called blocks of the graph. Therefore  a separable graph must contain blocks. Separable graphs have useful properties to highly accelerate the computation of combinatorial problems by decomposing the entire graph into blocks, as can be seen in \cite{alfredo} and \cite{dn}. The following property clearly accelerates the computation of budget-constrained cutting problems in separable graphs.

\begin{prope} \label{bloque} The Budget-Constrained Min-Cut (Max-Cut) in a separable graph is the minimum (maximum) budget-constrained min-cut (max-cut) of each one of its blocks.
\end{prope}

\begin{propo} \label{maxBound}
Let $r_i=\frac{w_i}{c_i}$ be the weight-cost ratio for each edge $i \in e$. Let $r_{(1)}, ..., r_{(m)}$ be the $m$ ratios of $G$ sorted in non-increasing sequence. Let $a$ be the largest index such that $\sum_{i=1}^a r_{(i)} \leq T$. An upper bound for the Budget-Constrained Max-Cut Problem is $\sum_{i=1}^a w_{(i)} + (T-\sum_{i=1}^a c_{(i)} )\frac{w_{(a+1)}}{c_{(a+1)}}$.
\dem Since the weight-cost ratio represents the profit of cutting each edge, this approach is similar to the procedure for solving a continuous Binary Knapsack Problem. Of course, the cut of $a$ or $a+1$ edges, could not be a cut but this would imply that the proposed value $\sum_{i=1}^a w_{(i)} + (T-\sum_{i=1}^a c_{(i)} )\frac{w_{(a+1)}}{c_{(a+1)}}$ is an upper bound. $\vspace{0.3cm}\hfill \square$
\end{propo}

\begin{coro}
\label{minBound}
Let $r_i=\frac{w_i}{c_i}$ be the weight-cost ratio for each edge $i \in e$. Let $r_{(1)}, ..., r_{(m)}$ be the $m$ ratios of $G$ sorted in non-decreasing  sequence. Let $a$ be the largest index such that $\sum_{i=1}^a r_{(i)} \leq T$. A lower bound for the Budget-Constrained Min-Cut Problem is $\sum_{i=1}^a w_{(i)} + (T-\sum_{i=1}^a c_{(i)} )\frac{w_{(a+1)}}{c_{(a+1)}}$.
\end{coro}

\begin{observation}The optimal value of the Budget-Constrained Min s-t Cut does not coincide with the optimal value of the integral Budget-Constrained Max s-t Flow.
\end{observation}
\dem A  well-known property (see Ford and Fulkerson \cite{fordmaximal}) establishes the relation between the min s-t cut of a network and the max s-t flow.  Here, we wonder whether the same or some modified property holds for the case with a budget constraint. As we will go on to show that property does not hold in the budget-constrained case. Adapted from \cite{flujo2} and \cite{flujo3}, we can formulate the integral Budget-Constrained Max s-t Flow for undirected graphs as:

\begin{alignat}{2} 
\mbox{} & \max\ F_{st} \label{flujoMAX}\\
\mbox{\bf s.t.}\ & & \nonumber \\
& \displaystyle \sum_{(s,j) \in E} f_{sj}  = F_{st} & \kern 10pt    \label{f1}\\
& \displaystyle \sum_{(j,t) \in E} f_{jt}  = F_{st} & \kern 10pt    \label{f2}\\
& \displaystyle \sum_{(i,j) \in E} f_{ij} - \sum_{(i,j) \in E} f_{ji} = 0 & \kern 10pt i \in V \setminus\{s,t\}   \label{f3}\\
& \displaystyle f_{ij} \leq w_{ij} x_{ij} & \kern 10pt   (i,j) \in E \label{f4}\\
& \displaystyle f_{ji} \leq w_{ij} x_{ij} & \kern 10pt   (i,j) \in E \label{f5}\\
& \sum_{(i,j) \in E} c_{ij} x_{ij} \leq T \label{budgetF} \\
& x_{ij} \in \{0,1\} & (i,j) \in E \label{binF}
\end{alignat}

$F_{st}$ represents the max s-t flow whilst $f_{ij}$ the quantity of flow that circulates through each edge $(i,j) \in E$. On the other hand, $x_{ij}$ takes value $1$ if some flow unit circulates through edge $(i,j)$, and $0$ otherwise. The objective function and the constraints \eqref{f1} - \eqref{f5} are the usual of the simple max s-t flow model but \eqref{f4} and \eqref{f5} have been modified inserting variables $x_{ij}$ to detect when an edge is used. Finally, constraint \eqref{budgetF} limits the total cost and as can be seen in \eqref{binF}, the integrality of the variables $x_{ij}$ is required.\\

Since it is compulsory to pay even for the flow sent through arcs that do not belong to the set of arcs with maximum flow, which are those that correspond to the minimum cut, then the optimal value of model \eqref{flujoMAX} does not coincide with the optimal value of model \eqref{fobj}. \\

\begin{figure}[H]
	\begin{center}
		\includegraphics[width=0.40\textwidth]{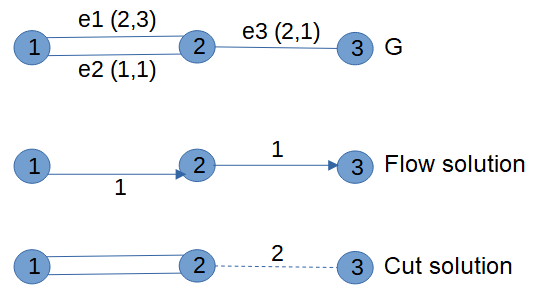}
	\end{center}
	\caption{Counterexample for budget-constrained flow and cut problems} 	\label{contraejemplo}
\end{figure}

We illustrate the above fact with a simple counterexample. It can be clearly seen in figure \ref{contraejemplo} for the graph shown at the top of this figure. The graph has three edges each one of them associated with one pair of numbers, the first one its weight and the second one its cost. In this example, we assume that the available budget is $\rho=2$. The solution for the budget-constrained max flow is shown in the central subfigure and it has an optimal value of $1$ with a cost of $2$ but the solution for the budget-constrained min cut, which is shown in the picture at the bottom of figure \ref{contraejemplo},  has an optimal value of $2$ with a cost of $1$. 

\section{Branch and bound approaches}

A first solution approach consists in solving the models \MIN\ and \MAX\ by integer programming. Besides using mathematical programming tools, in this section we propose an exact combinatorial algorithm to solve the budget-constrained min-and-max cut problems based on a branch and bound strategy.

The rationale behind this algorithm, while simple, is effective. It is based on the idea of introducing feasible cuts and then forbidding these cuts and introduce other feasible cuts in order to find and to quantify the objective value of all the feasible cuts. To do this, a branch and bound algorithm is developed and  the initial idea is also improved by adding the estimation of bounds to avoid exploring certain branches. Hence, this section shows an exact combinatorial algorithm for getting the best budget-constrained min and max cuts of a graph. For this, we look for feasible combinations of cuts bounded by an upper cost. The proposed method to obtain these combinations is explained in Subsection \ref{feasible_cut}. Moreover, the  exact and heuristic algorithms are described in Subsection \ref{exact_algorithm}.

\subsection{Finding a feasible cut for a budget constraint} \label{feasible_cut}

We have adapted the min-cut Stoer-Wagner algorithm \cite{mincut} in order to get any cut of a graph under a budget. The Stoer-Wagner algorithm splits the vertices of a graph into two clusters, $A$ and $B$. Initially, only an arbitrary node $a$ belongs to $A$ but iteratively $A$ grows by adding the most tightly connected node at each step. When $A=V$, this procedure called $MininumCutPhase$ returns the  minimum $s-t$ cut of the last two vertices $s$ and $t$ added to $A$. So, the main procedure called $MinimumCut$ goes on merging these two vertices and calling $MinimumCutPhase$ again. Once all the nodes are merged after $|V|$ iterations, the minimum cut found is the minimum cut of the graph.\\

At each step, the Stoer-Wagner algorithm evaluates the cut that splits the nodes of cluster $A$ from $B$ and if this cut is better than the best one saved, this is updated as the incumbent minimum cut. Instead, our adaptation consists in stopping the algorithm once a feasible cut is found, since we are not looking for a minimum cut but for any budget one. Obviously, if there is no feasible cut, as many iterations as Stoer-Wagner algorithm will be run and no cut will be returned, since even a minimum cut is not feasible. This implies that the computational complexity of the modified algorithm is the the same as the computational complexity of the Stoer-Wagner algorithm, $\mathcal{O}(|E| + |V| log |V|)$, since this can happen in the worst case. We have called the adapted procedures $BudgetedCutPhase$ and $BudgetedCut$ and the changes are highlighted in grey in Algorithms \ref{st1} and \ref{st2}. $\rho$ is the admissible budget and $cost(cut(A))$ means the total cost of dropping the edges of the cut.\\
 
\begin{algorithm}[h]
	\setcounter{AlgoLine}{0}
	\label{st1}
	\caption{$BudgetedCutPhase(G, a,\rho)$}
	$A \leftarrow \{a\}$ \\
	\While{$A \neq V$ \textcolor{arsenic} {\textbf{and} $cost(cut(A)) > \rho$} }
	{
		add to $A$ the most tightly connected vertex\\
	}
	shrink $G$ by merging the two vertices added in the last step\\
	\textbf{return} the current cut-of-the-phase\\
\end{algorithm}

\begin{algorithm}[h]
	\setcounter{AlgoLine}{0}
	\label{st2}
	\caption{$BudgetedCut(G, \rho)$}
	\While{$|V|>1$ \textbf{and} not budgeted cut is found}
	{
		$BudgetedCutPhase(G, v_1,\rho)$\\
	}
	\textbf{return} the budgeted cut if this is found, $\emptyset$ on the contrary
\end{algorithm}

\begin{figure}[H]
	\begin{center}
		\includegraphics[width=0.4\textwidth]{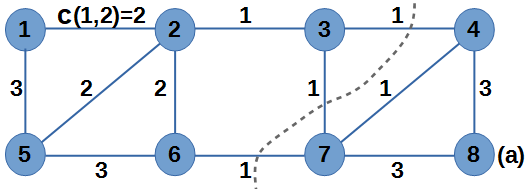}
	\end{center}
	\caption{Adaptation of the Stoer-Wagner algorithm} 	\label{ejemplo0}
\end{figure}

Figure \ref{ejemplo0} shows how the $BudgetedCut$ procedure works versus $MinimumCut$. The numbers on the edges are the costs of every cut-edge, their weights being irrelevant for illustrating the process. Let $\rho=3$ be the budget and $a=8$ the initial node for both procedures. When $A=\{8\}$ the cost of the cut is $6$ ($cost(cut(A))=6$) which is not within the budget. Both procedures iteratively grow the cluster $A$ for the first two steps as follows $A=\{8, 4\}$ ($cost(cut(A))=5$) and $A=\{8, 4, 7\}$ ($cost(cut(A))=3$). It is clear that at this step a cut under the budget is found, the cut that separates $A=\{8, 4, 7\}$ from the remaining nodes. Note that it was obtained by the $BudgetedCutPhase$ and returned by the $BudgetedCut$ procedure and its feasibility implies the end of the adapted algorithm, although the Stoer-Wagner algorithm would still have called $MinimumCutPhase$ $n-1=7$ times because it is not sure when a found cut is minimal.\\

\subsection{Branch and bound algorithm} \label{exact_algorithm}

\begin{figure}[H]
	\begin{center}
		\includegraphics[width=1\textwidth]{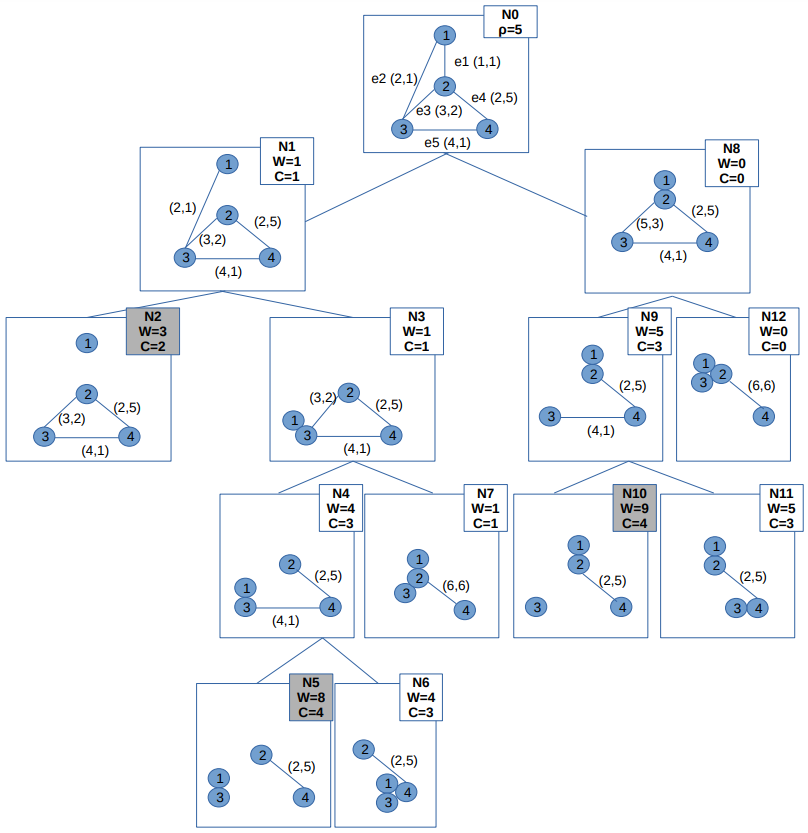}
	\end{center}
	\caption{Branch and bound example for min-cut} 	\label{ejemploMin}
\end{figure}

In this subsection, we describe the branch and bound algorithm and we propose to find an exact solution to the minimum and the maximum cuts of a graph bounded by a budget. The main idea is to get any cut within a budget then compute its goodness and, after banning the already evaluated cut, go on to get new budgeted cuts. Banning a cut implies performing a non-polynomial action since the way to ban it is to run all the combinations of allowed-forbidden cuts for every edge belonging to the found cut. Figure \ref{ejemploMin} illustrates how the method works, in an intuitive way. The graph on which it has been applied is show at the top of the figure, this is, we have a graph of $4$ vertices and $5$ edges. The weight and the cost of cutting each edge are respectively represented by the first and the second element of the pair into brackets, so the weight of edge $e2$ is $2$ and the cost for cutting it is $1$. The available budget is $\rho=5$. Each picture represents one node of the branch and bound tree and its box contains information about the order in which the node is explored (the label with the number of the node), the total cut weight (W) and the cost of this solution (C). Besides, when the represented solution is feasible, the box is shaded.\\

So, our method begins searching for a cut under the budget of $5$. The procedure $BudgetedCut$ is called and it returns the cut $\{e_1, e_2\}$, then two branches are sequentially added to the tree numbered N1 and N2. The cut of edge $e_1$ originates $N1$, the first node solution tree and $N2$ is originated by adding the cut of edge $e_2$. Note how $N1$ is not feasible which means that it is still not a cut but the cut of $e_1$ carries a weight of $1$ and a cost of $1$. Regarding $N2$, the found cut $\{e_1, e_2\}$ is fully done so it represents a feasible solution with a value $W=3$ and a cost $C=2$. Since we still do not know if this is optimal, the algorithm must continue until exploring all the possible solutions. Then, in following steps, both cuts will be alternatively and successively banned in order to get more feasible solutions from all the possible combinations. So, we reverse the last cut edge $e_2$ by banning its cut.  If we represent by $\overline{e}$ the cut of an edge $e$ and by $\underline{e}$ its contraction, we have already explored the combinations $N1=\{\overline{(1,2)}\}$, $N2=\{\overline{(1,2)}, \overline{(1,3)} \}$ and $N3=\{\overline{(1,2)}, \underline{(1,3)} \}$. As is usual, banning an edge means its two nodes are merged. Thus, banning the cutting of $(1,3)$ originates solution $N3$  which is non-feasible since it does not cut the graph. Then from $N3$ we try to get another cut under the remaining budget $4$. If necessary, weights and costs of edges are increased when nodes are shrunk as it happens in nodes $N7$, $N8$ and $N12$. The cut $\{e_3, e_5\}$ is returned when $BudgetedCut$ Procedure is called from $N3$ and the algorithm follows the same path.\\

In addition to Figure \ref{ejemploMin}, Table \ref{tEjemploMin} also illustrates the example showing the summary of all the branch and cut combinations explored. On the other hand, the process for searching new cuts each time a cut is fully applied, consists in discarding the cut by successively shrinking the nodes of each of its edges and thus finding another cut. Table \ref{tEjemploMin} shows the nodes for which $BudgetedCut$ Procedure is called and its returned cuts. The columns $Weight$ and $Cost$, respectively, are the weight and the cost of each new cut found without summing the weight and cost of the edges already cut.\\

\begingroup
\setlength{\tabcolsep}{10pt}
\renewcommand{\arraystretch}{1.5}
\begin{table}[H]
\centering
\begin{tabular}{|l|l|l|r|r|}
\hline
\textbf{Node} &\textbf{Parent} &\textbf{Cuts / Shrinks} & \textbf{Weight} & \textbf{Cost}  \\ \hline
$N1$ & $N0$ & $\overline{(1,2)}$ & $1$ & $1$ \\ \hline
\cellcolor{Gray}$N2$ & \cellcolor{Gray} $N1$ & \cellcolor{Gray} $\overline{(1,2)}$,  $\overline{(1,3)}$ & \cellcolor{Gray} $3$ & \cellcolor{Gray} $2$ \\ \hline
$N3$ & $N1$ & $\overline{(1,2)}$, $\underline{(1,3)}$ & $1$ & $1$ \\ \hline
$N4$ & $N3$ & $\overline{(1,2)}$, $\underline{(1,3)}$, $\overline{(3,2)}$ & $4$ & $3$ \\ \hline
\cellcolor{Gray}$N5$ & \cellcolor{Gray} $N4$& \cellcolor{Gray}$\overline{(1,2)}$,$\underline{(1,3)}$, \cellcolor{Gray} $\overline{(3,2)}$, $\overline{(3,4)}$ & \cellcolor{Gray} $8$ & \cellcolor{Gray} $4$ \\ \hline
$N6$ & $N4$& $\overline{(1,2)}$, \underline{$(1,3)$}, $\overline{(3,2)}$, \underline{$(3,4)$} & $4$ & $3$ \\ \hline
$N7$ & $N3$ & $\overline{(1,2)}$, \underline{$(1,3)$}, \underline{$(3,2)$} & $1$ & $1$ \\ \hline
$N8$ & $N0$ & \underline{$(1,2)$} & $0$ & $0$ \\ \hline
$N9$ & $N8$ & \underline{$(1,2)$}, $\overline{(2,3)}$ & $5$ & $3$ \\ \hline
\cellcolor{Gray}$N10$ & \cellcolor{Gray}$N9$ & \cellcolor{Gray}\underline{$(1,2)$}, $\overline{(2,3)}$, $\overline{(3,4)}$ & \cellcolor{Gray}$9$ & \cellcolor{Gray}$4$ \\ \hline
$N11$ & $N9$ & \underline{$(1,2)$}, $\overline{(2,3)}$, \underline{$(3,4)$} & $5$ & $3$ \\ \hline
$N12$ & $N0$ & \underline{$(1,2)$}, \underline{$(2,3)$} & $0$ & $0$ \\ \hline
\end{tabular}
\caption{Nodes of the exact algorithm example}
\label{tEjemploMin}
\end{table}
\endgroup

\begingroup
\setlength{\tabcolsep}{10pt}
\renewcommand{\arraystretch}{1.5}
\begin{table}[H]
\centering
\begin{tabular}{|l|l|l|r|r|}
\hline
\textbf{Node} & \textbf{Already cut} & \textbf{Cut} & \textbf{Weight} & \textbf{Cost}  \\ \hline
$N0$ & $\emptyset$ & $(1,2)$, $(1,3)$ & $4$ & $3$ \\ \hline
$N3$ & $(1,2)$ & $(2,3)$, $(3,4)$ & $7$ & $3$ \\ \hline
$N6$ & $(1,2)$ & $\emptyset$ & $0$ & $0$ \\ \hline
$N7$ & $(1,2)$ & $\emptyset$ & $0$ & $0$ \\ \hline
$N8$ & $\emptyset$ & $(2,3)$, $(3,4)$ & $9$ & $4$ \\ \hline
$N11$ & $(3,4)$ & $\emptyset$ & $0$ & $0$ \\ \hline
$N12$ & $\emptyset$ & $\emptyset$ & $0$ & $0$ \\ \hline
\end{tabular}
\caption{All the cuts returned by $BudgetedCut$ Procedure}
\label{tEjemploCuts}
\end{table}
\endgroup

As it can be seen, the best solution is $3$. Its optimality is certified since all the possible cuts were explored. Later, we shall improve the performance of the algorithm by introducing bounds that allow to implicitly enumerate some cuts which will never be optimal. We also note that in the above example the optimal value would have been achieved  by just calling the Stoer-Wagner Algorithm since the cost of the solution of minimum weight is directly under the budget. However, we presented the full exploration  only for illustrating the proposed method.\\

\begin{figure}[H]
	\begin{center}
		\includegraphics[width=1\textwidth]{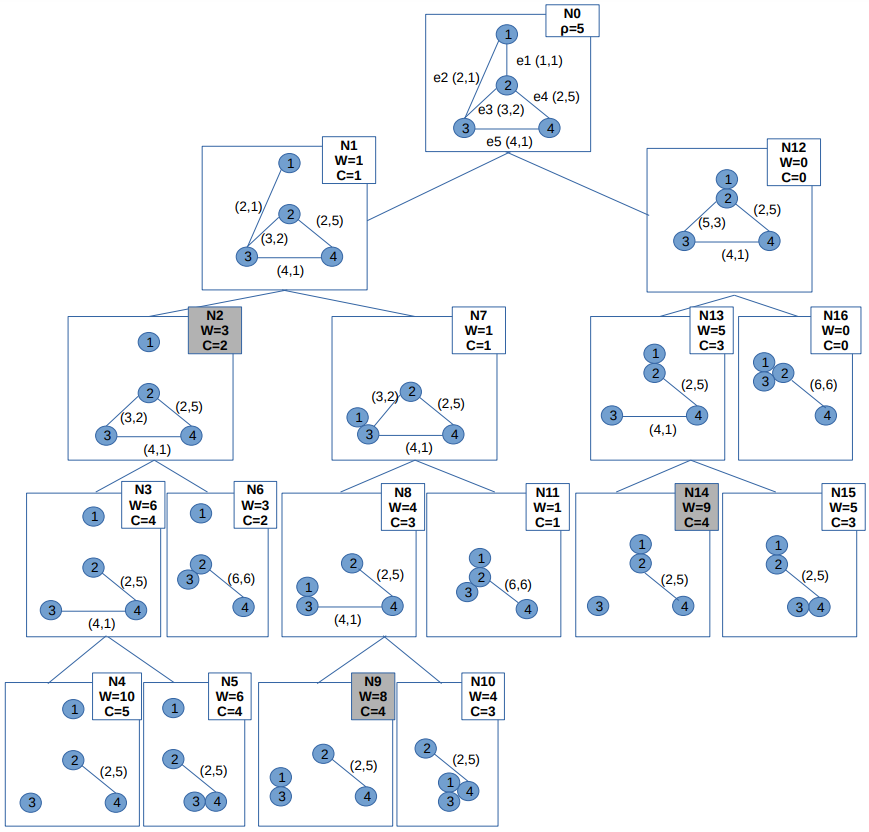}
	\end{center}
	\caption{Branch and bound example for max-cut} 	\label{ejemploMax}
\end{figure}

Figure \ref{ejemploMax} illustrates the method for the max cut problem. Basically, it is identical except for the fact that once a cut is found we have to go on searching for more cuts instead of going backward into the tree, in order to maximize the total weight of the cut. So, $N3, N4, N5$ and $N6$ of the branch and bound tree of the max cut problem are not explored in the minimum cut case.\\

Besides Figure \ref{ejemploMax}, Table \ref{tEjemploMax} also illustrates the example showing the summary of all the branch and cut combinations explored for the max cut problem. 

\begingroup
\setlength{\tabcolsep}{10pt}
\renewcommand{\arraystretch}{1.5}
\begin{table}[H]
\centering
\begin{tabular}{|l|l|l|r|r|}
\hline
\textbf{Node} &\textbf{Parent} &\textbf{Cuts / Shrinks} & \textbf{Weight} & \textbf{Cost}  \\ \hline
$N1$ & $N0$ & $\overline{(1,2)}$ & $1$ & $1$ \\ \hline
\cellcolor{Gray}$N2$ & \cellcolor{Gray} $N1$ & \cellcolor{Gray}$\overline{(1,2)}$, $\overline{(1,3)}$ & \cellcolor{Gray} $3$ & \cellcolor{Gray} $2$ \\ \hline
$N3$ & $N2$ & $\overline{(1,2)}$, $\overline{(2,3)}$ & $6$ & $4$ \\ \hline
$N4$ & $N3$ & $\overline{(1,2)}$, $\overline{(2,3)}$, $\overline{(3,4)}$ & $10$ & $5$ \\ \hline
$N5$ & $N3$ & $\overline{(1,2)}$, $\overline{(2,3)}$, \underline{$(3,4)$} & $6$ & $4$ \\ \hline
$N6$ & $N2$ & $\overline{(1,2)}$, \underline{$(2,3)$}, $\overline{(1,3)}$ & $3$ & $2$ \\ \hline
$N7$ & $N1$ & $\overline{(1,2)}$, \underline{$(1,3)$} & $1$ & $1$ \\ \hline
$N8$ & $N7$ & $\overline{(1,2)}$, \underline{$(1,3)$}, $\overline{(3,2)}$ & $4$ & $3$ \\ \hline
\cellcolor{Gray}$N9$ & \cellcolor{Gray} $N8$& \cellcolor{Gray}$(1,2)$, \cellcolor{Gray} \underline{$(1,3)$}, \cellcolor{Gray} $\overline{(3,2)}$, $\overline{(3,4)}$ & \cellcolor{Gray} $8$ & \cellcolor{Gray} $4$ \\ \hline
$N10$ & $N8$& $\overline{(1,2)}$, \underline{$(1,3)$}, $\overline{(3,2)}$, \underline{$(3,4)$} & $4$ & $3$ \\ \hline
$N11$ & $N7$ & $\overline{(1,2)}$, \underline{$(1,3)$}, \underline{$(3,2)$} & $1$ & $1$ \\ \hline
$N12$ & $N0$ & \underline{$(1,2)$} & $0$ & $0$ \\ \hline
$N13$ & $N12$ & \underline{$(1,2)$}, $\overline{(2,3)}$ & $5$ & $3$ \\ \hline
\cellcolor{Gray}$N14$ & \cellcolor{Gray}$N13$ & \cellcolor{Gray}\underline{$(1,2)$}, $\overline{(2,3)}$, $\overline{(3,4)}$ & \cellcolor{Gray}$9$ & \cellcolor{Gray}$4$ \\ \hline
$N15$ & $N13$ & \underline{$(1,2)$}, $\overline{(2,3)}$, \underline{$(3,4)$} & $5$ & $3$ \\ \hline
$N16$ & $N0$ & \underline{$(1,2)$}, \underline{$(2,3)$} & $0$ & $0$ \\ \hline
\end{tabular}
\caption{Nodes of the exact algorithm example}
\label{tEjemploMax}
\end{table}
\endgroup
Table \ref{tEjemploMaxCuts} show the nodes in which the $BudgetedCut$ Procedure is called and its returned cuts.

\begingroup
\setlength{\tabcolsep}{10pt}
\renewcommand{\arraystretch}{1.5}
\begin{table}[H]
\centering
\begin{tabular}{|l|l|l|r|r|}
\hline
\textbf{Node} & \textbf{Already cut} & \textbf{Cut} & \textbf{Weight} & \textbf{Cost}  \\ \hline
$N0$ & $\emptyset$ & $(1,2)$, $(1,3)$ & $4$ & $3$ \\ \hline
$N2$ & $(1,2)$, $(1,3)$ & $(2,3)$, $(3,4)$ & $7$ & $3$ \\ \hline
$N4$ & $(1,2)$, $(1,3)$, $(2,3)$, $(3,4)$ & $\emptyset$ & $0$ & $0$ \\ \hline
$N7$ & $(1,2)$ & $(2,3)$, $(3,4)$ & $7$ & $3$ \\ \hline
$N10$ & $(1,2)$ & $\emptyset$ & $0$ & $0$ \\ \hline
$N11$ & $(1,2)$ & $\emptyset$ & $0$ & $0$ \\ \hline
$N12$ & $\emptyset$ & $(2,3)$, $(3,4)$ & $9$ & $4$ \\ \hline
$N15$ & $(3,4)$ & $\emptyset$ & $0$ & $0$ \\ \hline
$N16$ & $\emptyset$ & $\emptyset$ & $0$ & $0$ \\ \hline
\end{tabular}
\caption{All the cuts found in the Max Cut Problem}
\label{tEjemploMaxCuts}
\end{table}
\endgroup

\begin{figure}[H]
	\begin{center}
		\includegraphics[width=0.75\textwidth]{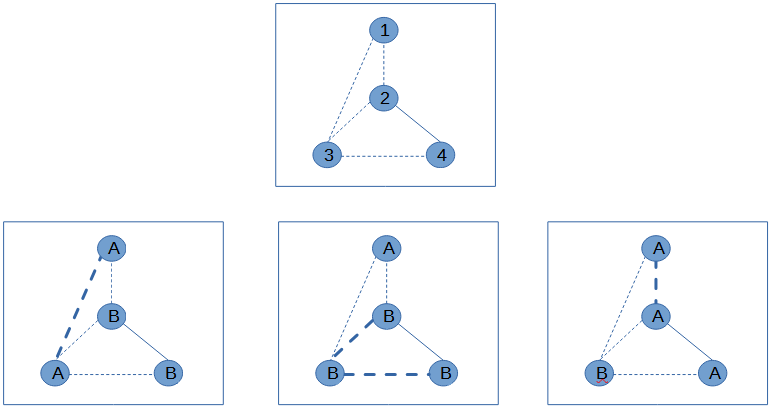}
	\end{center}
	\caption{Non-valid cut} 	\label{noValido}
\end{figure}

It must be noted that for this case $N4$ is not a solution despite coming from a new cut. That is analyzed by Figure \ref{noValido}. The picture at the top of the figure represents the graph. The cut edges have been represented by dotted lines. The three pictures at the bottom of the figure represent all the combinations of the nodes divided between two clusters $A$ and $B$. Nodes $2$ and $4$ must belong to the same cluster since both are connected. Note how for any possible combination there is at least one cut edge between two nodes belonging to the same cluster that has been represented by a thick dotted line. So, before considering the new cut its feasibility needs to be checked.  \\

The final procedure is described by Algorithm \ref{main}, which also calls the forward procedure (Algorithm \ref{forward}) and the backward procedure (Algorithm \ref{backward}). The design of Algorithm \ref{main} is the typical of a depth-first search tree and it basically alternates between a forward phase in which new branches are added to the tree and a backward phase in which a branch is reversed and others may also be deleted in order to explore different solutions. That strategy optimizes the required computational memory since the solution space could be too vast since it is only required to store the last solution, in the successive also called current solution, as well as the cuts and merges that lead to it. So, as long as not all the possible solutions are explored, a bound from the current solution is estimated in each loop. In the case of the minimum case, just by getting a minimum cut without budget constraint since that has negligible computational cost. Nevertheless in the maximum case getting a similar bound (not constrained) would be NP-Hard. Therefore, we estimate it from Proposition \ref{maxBound}. If the bound is better than the best solution being, lower than the best one in the minimum case but higher in the maximum case, then the algorithm goes forward by searching a new cut. Otherwise, the algorithm goes backward by reversing the last cut edge. Note also that in the maximum case the algorithm goes on searching for new cuts before calling the Backward Procedure as long as a feasible cut is found. \\

Remember, not every new cut satisfying the budget is valid as was illustrated in Figure \ref{noValido}. In practice, checking the feasibility of a new cut can be easily done and besides it is computationally cheap by selecting any node and labelling it as belonging to one of the clusters. We go on labelling nodes by DFS (Depth-First Search) or BFS (Breadth-First Search), for both cut and uncut edges. If a neighboring node belongs to an uncut edge, then it has to be labelled as belonging to the same cluster, otherwise to the other cluster. If all the nodes are labelled without detecting inconsistencies about the cluster to which they have to belong, then the cut is valid otherwise, the cut is not valid.\\

\begin{algorithm}[H] \label{main}
	\setcounter{AlgoLine}{0}
	\caption{$ConstrainedCut$}
	\SetKwInOut{Input}{Input parameters}
	\SetKwInOut{Output}{Output}
	\Input{\\
	$G$: the graph to be partitioned,\\
	$\rho$: budget,\\
	$type:$ the problem type ($MIN$, $MAX$)
     \rule{0.8\paperwidth}{0.4pt}
	}
    \vspace{0.25cm}
        \textbf{if} $type=MIN$ \textbf{then} $best=+\infty$ \textbf{else} $best=0$\\
        $R=\rho$ \tcp*{$R$: the available budget for new cuts}
	\While{the solutions tree is not explored}
	{
            \textbf{if} $bound$ better than $best$ \textbf{ then } $X=BudgetedCut(G, R)$; \textbf{ else } $X=\emptyset$;\\
            \While{$X \neq \emptyset$}
	    {                
                \textbf{if} the current solution is valid and better than $best$ \textbf{then} save it and update $best$;\\
                $forward(X)$; \\
                \textbf{if} $type=MIN$ or it is not a valid solution \textbf{ then } $X=\emptyset$; \textbf{ else } $X=BoundedCut(G, R)$;\\
            }
            $backward()$;
	}
\end{algorithm}

Regarding the $Forward$ Procedure, it just adds a branch to explore for each edge belonging to the cut returned by the $BoundedCut$ Procedure and it decreases the available budget for future cuts. As to $ForwardCut$ Procedure it removes the deepest node of the tree if this has been forbidden or shrunk until a node corresponding to an edge with cut state is found. Then, the cut of this edge is reversed by shrinking it.\\

\begin{algorithm}[H] \label{forward}
	\setcounter{AlgoLine}{0}
	\caption{$forward(X)$}
        \tcp{the cut is applied}
        \ForEach{$e \in X$}
        {
            $G = G \setminus e$  \\
            add $e$ with state '$cut$' to the tree\\
            $R=R-cost(e)$
        }
\end{algorithm}

\begin{algorithm}[H] \label{backward}
	\setcounter{AlgoLine}{0}
	\caption{$backward()$}
	        Let $e=(u,v)$ be the edge of the deepest node from the tree\\
    	    \While{state of $e$ $=$ '$shrunken$'}
    	    {
    	        remove the deepest node from the tree\\
    	        restore $e$ in the graph\\
    	    }
    	    \If{$\exists e$}
    	    {
    	        \tcp{the tree is not empty then the state of $e$ is '$cut$'}
    	        change the state of $e$ to '$shrunken'$\\
    	        shrink $u$ and $v$ in the graph\\
    	    } 
\end{algorithm}

\section{Lagrangian relaxation approaches}

Lagrangian relaxation allows to provide tight lower bounds for the  the budgeted min-cut problem, solving a pseudopolynomial number of  simpler min-cut problems each with polynomial complexity. For considering the Lagrangian relaxation, it is enough to relax constraint (\ref{budget}) by defining the Lagrangian function as {$L_x(\lambda) = l_x^1 + \lambda \hspace{0.1cm} l_x^2$ where $l_x^1=\sum_{(i,j) \in E} w_{ij} x_{ij}$ and $l_x^2=\sum_{(i,j) \in E} c_{ij} x_{ij} - T$}. $L_x(\lambda)$ will be the new objective function of the successive subproblems to solve by fixing the multiplier {$\lambda\ge 0$}. The corresponding Lagrangian dual $\max_{\lambda\ge 0} \min_x L_x(\lambda)$ is always a very good lower bound and in most cases results in the optimal solution of the budged constrained min-cut problem.  Regarding the budget-constrained max-cut problem, removing constraint (\ref{budgetMax}) does not simplify the problem which remains NP-Complete. On the contrary, strengthening the max-cut problem reduces its solution space which usually implies less computational effort. So, we will analyse the use of the Lagrangian relaxation only for Budget-Constrained Min-Cut Problem.\\

The first approach to solve the Lagrangian dual solves each Lagrangian subproblem, with objective function {$L_x(\lambda)$} for each $\lambda>0$, using an integer programming solver. Algorithm \ref{lagrangian} explains the strategy followed in order to iterate on the different values of $\lambda$ considered for the  the different relaxed subproblems.

\begin{algorithm}[H] \label{lagrangian}
	\setcounter{AlgoLine}{0}
	\caption{Lagrangian approach}
	$\lambda_A=0$\\
        $A$= solution of relaxed \MIN\ with $\lambda_A$\\
        \If{$A$ is feasible}
        {
            $optimal=A$\\
        }
        \Else{
            $\lambda_A=0$\\                     
            $\lambda_B=$ big value\\
            $B$= solution of relaxed \MIN\ with $\lambda_B$\\
            \While{new solutions are found}
            {            
                $\displaystyle \lambda=\frac{l_B^1 - l_A^1}{l_A^2 - l_B^2}$\\
                $C$= solution of relaxed \MIN\ with $\lambda$\\                
                \If{$C$ is feasible}{
                    $\lambda_B=\lambda$
                }
                \Else{
                    $\lambda_A=\lambda$
                }
            }
        }
\end{algorithm}

We have followed the usual one-dimensional search algorithm for problems with only one relaxed constraint. Then, we hold two solutions at each iteration. Solution $A$ is a non-feasible solution in the sense that it is a cut exceeding the budget. Solution $B$ is any cut that does not exceed the budget but it is non necessarily the one giving the value of the Lagrangian dual. $A$ and $B$ will be progressively adjusted throughout the iterations. So, initially $A$ is the solution reached with $\lambda_A=0$, that is, without budget restriction,  and the initial $B$ solution will be obtained by fixing a big value to $\lambda_B$, that guarantees a feasible solution, although most likely it will not be optimal. Thus, we should find another multiplier $\lambda$ trying to get an optimal solution for the dual. To this, as is usual,  we fix the new multiplier according to $\displaystyle \lambda=\frac{l_B^1 - l_A^1}{l_A^2 - l_B^2}$ to search new solutions in an orthogonal direction to $\overline{AB}$. If we get a new solution $C$ then we fix $\lambda_B=\lambda$ if $C$ is feasible but $\lambda_A=\lambda$ if $C$ is over the budget. Hence, in each iteration we are reducing the search space of the multiplier. The algorithm ends when no new solution is found. Note we only identify solutions on the convex hull of the \MIN, so the optimality is not assured because there may be solutions in the non-convex space, so this algorithm is not exact.\\

Given that removing constraint \eqref{budget} the new problem is just a simple min-cut problem, one can also get its solution by the Stoer-Wagner algorithm. {Furthermore, in this way, other cuts are obtained during the process of getting each min-cut. Any of these non-final min-cuts of the Lagrangian subproblem may be feasible and even better for the original problem than the final min-cut. Then, by evaluating them, the solution can be improved and in some instances we will get better solutions than those returned by the IP Lagrangian relaxation approach}. This is our second method to solve Lagrangian relaxation that we analyze in the computational experiment section.

\section{Computational results}

By this computational experience we assess and compare the efficiency of the methods proposed to solve both (min-and-max) budget-constrained cuts. With that aim, we have used random connected graphs with $n$ vertices and $m$ edges, already used in \cite{Miyazawa_etal_2021} \footnote{https://www.loco.ic.unicamp.br/files/instances/bcpk/}. Table \ref{instancias} shows the sizes and number of these instances used in our experiments.\\

\begin{table}[h]
\centering
\begin{tabular}{|r|r|r|}
\hline
\textbf{n} & \textbf{m} & \textbf{\# instances}\\
\hline
\rowcolor{Gray}
\textbf{20} & & \textbf{90} \\
 & 30 & 30 \\
 & 50 & 30 \\
 & 100 & 30 \\
\rowcolor{Gray}
\textbf{30} & & \textbf{90} \\
 & 50 & 30 \\
 & 70 & 30 \\
 & 200 & 30 \\
\rowcolor{Gray}
\textbf{50} & & \textbf{90} \\
 & 70 & 30 \\
 & 100 & 30 \\
 & 400 & 30 \\
 \rowcolor{Gray}
\textbf{70} & & \textbf{90} \\
 & 100 & 30 \\
 & 200 & 30 \\
 & 600 & 30 \\
 \rowcolor{Gray}
\textbf{100} & & \textbf{90} \\
 & 150 & 30 \\
 & 300 & 30 \\
 & 800 & 30 \\
 \rowcolor{Gray}
\textbf{200} & & \textbf{90} \\
 & 300 & 30 \\
 & 600 & 30 \\
 & 1500 & 30 \\
 \hline
\end{tabular}
\caption{Size of instances}
\label{instancias}
\end{table}

For each instance, the weights and costs of each edge were randomly set from $1$ to $10$. Besides, once we compute the simple minimum cut according to the weights ($m_w$) of each instance and also the weight of the simple minimum cut according to the costs ($m_c$), different budgets were fixed under the expression $T=m_w+p$ $m_c$. Min-cut problems were solved for the values $p=25\%$, $50\%$ and $75\%$ whilst max-cut problems were solved for $p=25\%$, $50\%$, $75\%$, $100\%$ and $200\%$. Our experiments were conducted on a PC with a 2.30 GHz Intel Xeon, 64 GB of RAM and we used as optimization engine IBM ILOG CPLEX Interactive Optimiser 20.1.0.0, setting a max computing time per instance of 3600 s. In the following, we summarize several aspects of the computational experience divided into two sections one for the min version of the problems and the other for the max version.\\

\subsection{Budget-Constrained Min-Cut computational experience.}

Table \ref{Tmin} reports results on computing time of the IP formulation (\textit{t LP}) and the branch and bound algorithm (\textit{t Algorithm}). Since not all the instances were solved by the IP formulation after 3600 seconds of computation, the last column reports the number of instances solved. Times reported are the average among all the instances considered. The required processing time presents a large deviation between both methods: 769.22 seconds for the IP formulation vs. 1.22 seconds for the algorithm.\\

It is also important to point out some details of the integer programming computation. Since the Stoer-Wagner algorithm can provide a feasible solution of minimum cost practically without any computational effort, it was given as an initial solution  to the solver in order to accelerate the computation and also to avoid instances without feasible solutions. In addition, given that the proportion of violated constraints  \eqref{corte2} - \eqref{corte3} is rather low, they have been handled as lazy constraints, which highly accelerated the computation.\\

\begin{table}[h]
    \centering
    \begin{tabular}{|r|r|r|r|r|r|}
    \hline
        \textbf{\% p} & \textbf{$\#$ Nodes} & \textbf{t IP} & \textbf{t Algorithm}  & \textbf{\# Conc.} \\ \hline
        \rowcolor{Gray}
        25 & ~ & 135.46 & 0.08 & 506 \\ \hline        
        ~ & 20 & 0.00 & 0.00 & 90\\ \hline
        ~ & 30 & 0.04 & 0.00 & 90\\ \hline
        ~ & 50 & 0.79 & 0.02 & 90\\ \hline
        ~ & 70 & 5.28 & 0.02 & 90\\ \hline
        ~ & 100 & 35.90 & 0.09 & 90\\ \hline
        ~ & 200 & 1156.54 & 0.55 & 56\\ \hline
        \rowcolor{Gray}
        50 & ~ & 170.48 & 0.40 & 516 \\ \hline        
        ~ & 20 & 0.01 & 0.00  & 90\\ \hline
        ~ & 30 & 0.04 & 0.00 & 90 \\ \hline
        ~ & 50 & 0.79 & 0.02 & 90 \\ \hline
        ~ & 70 & 5.19 & 0.06  & 90\\ \hline
        ~ & 100 & 34.49 & 0.18 & 90 \\ \hline
        ~ & 200 & 1277.62 & 2.74 & 66\\ \hline        
        \rowcolor{Gray}
        75 & ~ & 159.73 & 1.15 & 516\\ \hline        
        ~ & 20 & 0.00 & 0.00 & 90\\ \hline
        ~ & 30 & 0.06 & 0.00 & 90\\ \hline
        ~ & 50 & 0.98 & 0.06 & 90\\ \hline
        ~ & 70 & 5.29 & 0.13 & 90\\ \hline
        ~ & 100 & 33.73 & 0.36 & 90\\ \hline
        ~ & 200 & 1194.17 & 8.27 & 66\\ \hline        
        \rowcolor{Gray}
        ~ & ~ & \textbf{769.22} & \textbf{1.22} & \textbf{1538} \\ \hline
    \end{tabular}
    \caption{Solved instances by both exact methods (min-cut)} \label{Tmin}
\end{table}

Table \ref{TNonmin} details the non-concluded instances. First of all, there are $82$ out of $1620$ instances for which the IP solver could not certify optimality after 3600 seconds, although 20 of them already achieved the optimum. However,  the exact branch and bound algorithm was able to conclude all the experiments with an average computing time  of $1.37$ seconds. Regarding the quality of the solutions of the non-concluded instances, Column $\% GAP$ quantifies the relative gap between the IP solution and the branch and bound algorithm according to $\% GAP=(v_{IP}-v_{Alg})*100/v_{Alg}$ where $v_{IP}$ is the objective value of the best solution returned by the solver and $v_{Alg}$ the value returned by the algorithm. The $\% GAP$ was $161.58 \%$ on average which means that, in general, the IP solution is more than twice the optimal solution for the non-concluded instances. Finally, although one initial feasible solution was provided to the solver, we also have results about the performance when no initial solution was supplied: the solver was not even capable of getting a feasible solution for $13$ out of the 20 non-concluded instances. (See Column \textit{\# Non-Solution}.)\\

\begin{table}[!ht]
    \centering
    \begin{tabular}{|r|r|r|r|r|r|r|}
    \hline
        \textbf{\% p} & \textbf{$\#$ Nodes} & \textbf{\# Non-concluded} & \textbf{t Algorithm} & \textbf{$\#$ Coincidents}   & \textbf{\% GAP}  
 & \textbf{\# Non-Solution} \\ \hline
        \rowcolor{Gray}
        ~ & \textbf{200} & \textbf{82} & \textbf{1.37} & \textbf{20}  & \textbf{161.58} & \textbf{13} \\ \hline        
        25 & 200 & 34 & 0.79 & 10   & 150.53 & 8 \\ \hline
        50 & 200 & 24 & 0.77 & 3   & 159.58 & 4 \\ \hline
        75 & 200 & 24 & 2.56 & 7   & 175.82 & 1 \\ \hline        
    \end{tabular}
    \caption{Non-solved instances by IP (min-cut)}  \label{TNonmin}
\end{table}

We made another attempt to accelerate the integer programming process trying to add valid knapsack cover inequalities from the Constraint \eqref{budget}. Those were added each time a feasible, but fractional, solution was found during the computation, but finally this does not improve the total computing time. A possible explanation for this is that few feasible fractional solutions were found. For example, solving the instance \textit{rnd\_20\_30\_a\_1} with $p=200 \%$ in the max-cut case, originated $54$ branch-and-bound branches of which only $3$ where feasible but fractional. For the remaining instances, the proportions were similar.\\

Regarding the Lagrangian relaxation Table \ref{TLang} summarizes the computational experience. Column \textit{t IP} shows the average time when the subproblems are solved by the IP solver, alternatively \textit{t S-W} reports results when the Stoer-Wagner Algorithm is used for getting the min-cut, \textit{\# Non-Opt} the number of returned non-optimal solutions and \textit{\% GAP} the relative gap between the optimal solutions and the solutions obtained by solving the Lagrangian relaxation subproblems by the Stoer-Wagner algorithm. The solutions obtained by the  Lagrangian relaxation were non-optimal only in $132$ out of $1620$ instances. Besides, the total average gap was $1.93 \%$.  On the other hand, we also have results about the performance when only the final min-cut is taken into account at each iteration. Then, $145$ non-optimal solutions and a total average gap of $2.03 \%$ are obtained. With regard to the computational time, all the instances were solved in $1$ second at most when the Stoer-Wagner algorithm was used to solve the subproblems. However, when the cuts were obtained by the IP solver, the time was even worse than when Lagrangian relaxation was not employed, so that no instance of $200$ nodes could be completed.  In this respect, it is important to comment certain details about the implementation of the Lagrangian relaxation when we use the IP routine. 

Recall that the constraints of \MIN\ problem, without the budget constraint, define the convex hull of the feasible domain. Thus, it is possible to solve the problem by linear programming. In spite of that, in our experiments, solving the linear model requires more time than solving the integer model, where the constraints (\ref{corte2}) are handled as lazy constraints. Certainly, this must be due to the large number (cubic) of constrains in the linear model. For this reason, \textit{t IP} reports the times solving the subproblems with the IP solver and lazy constraints.\\

 \begin{table}[!ht]
    \centering
    \begin{tabular}{|r|r|r|r|r|r|r|}
    \hline
        \textbf{\% p} & \textbf{$\#$ Nodes}  & \textbf{t IP} & \textbf{t S-W}  & \textbf{\# Non-Opt}  & \textbf{\% GAP}  \\ \hline
        \rowcolor{Gray}
        \textbf{25} &  & \textbf{614.10} & \textbf{0.02}   & \textbf{23} & \textbf{1.50} \\ \hline
          & 20 & 0.01 & 0.00   & 0 & 0.00\\ \hline
          & 30 & 0.08 &  0.00   & 3  & 0.41  \\ \hline
          & 50 & 1.03 & 0.00    & 7 & 1.42 \\ \hline
           & 70 & 7.6 & 0.01    & 5 & 3.16 \\ \hline         
          & 100 & 61.88 & 0.02    & 4 & 1.17 \\ \hline                  
           & 200  & 3600.00 & 0.12    & 4 & 2.87 \\ \hline   
         \rowcolor{Gray}
         \textbf{50} & & \textbf{614.11} & \textbf{0.03}    & \textbf{48} & \textbf{1.90} \\ \hline
          &  20 & 0.02 &   0.00    & 4 & 0.69 \\ \hline   
          &  30 &  0.09 & 0.00     & 3  & 0.34\\ \hline   
          &  50 & 1.20  & 0.00   & 9 & 2.04 \\ \hline             
          &  70  & 8.80 &   0.01   & 7 & 1.53 \\ \hline   
          &  100  & 74.55  & 0.02   & 13 & 3.25  \\ \hline             
          &  200  & 3600.00 & 0.14    & 12 &  3.57\\ \hline  
          \rowcolor{Gray}
        \textbf{75} & & \textbf{614.27} & \textbf{0.03}  & \textbf{61} & \textbf{2.40}  \\ \hline
          &  20  & 0.02 & 0.00  & 5  & 1.15\\ \hline   
          &  30  & 0.09 &  0.00   & 10 & 2.15 \\ \hline 
          &  50  & 1.03 & 0.01  & 9  & 2.08 \\ \hline           
          &  70  & 7.74 & 0.01  & 12 & 2.57 \\ \hline
          &  100  & 76.74 & 0.01  & 14 & 3.57 \\ \hline
          &  200 & 3600.00 & 0.17  & 11 & 2.83 \\ \hline
           \rowcolor{Gray}
          & & \textbf{614.16} & \textbf{0.03}   & \textbf{132} & \textbf{1.93}  \\ \hline
    \end{tabular}
    \caption{Lagrangian relaxation summary} \label{TLang}
\end{table}

An enhancement for improving the Lagrangian relaxation consists in implementing a local search procedure applied to each solution obtained. The inclusion of a basic local search, based on the Kernighan-Lin algorithm, leads to a modest refinement in results: $125$ non-optimal solutions discovered, accompanied by an average GAP percentage of $1.54 \%$ without worsening computational time.

\subsection{Budget-Constrained Max-Cut computational experience}

From an overall number of  $2700$ instances,  $1948$ instances were solved to optimality by both methods, among them $689$ were concluded by the algorithm but not for the IP solver, and $71$ were not solved to optimality by any of the methods. Table \ref{tMax} reports results on computing time of the IP formulation (\textit{t LP}) and the branch and bound algorithm (\textit{t Algorithm}). Since not all the instances were solved after 3600 seconds, only the instances solved by both methods are summarized, their number being reported in the last column. The cells in the table contain the average processing time in seconds. Solving the IP formulation required an average $207.30$ seconds vs. $65.25$ for the algorithm. Although the branch and bound algorithm clearly spent less computational time, the gaps are not as large as for the min-cut case. In addition, on this occasion, there are also instances not solved by the algorithm. Recall that even the ordinary max-cut problem is already NP-Complete.\\

\begin{table}[!ht]
    \centering
    \begin{tabular}{|r|r|r|r|r|r|r|}
    \hline
        \textbf{\% p} & \textbf{$\#$ Nodes} & \textbf{t LP} & \textbf{t Algorithm} & \textbf{\# Conc.} \\ \hline
        \rowcolor{Gray}
        25 & ~ & 208.10 & 0.01 & 422 \\ \hline
        ~ & 20 & 0.13 & 0.00 & 90 \\ \hline
        ~ & 30 & 1.67 & 0.01  & 90 \\ \hline
        ~ & 50 & 99.97 & 0.02  & 90 \\ \hline
        ~ & 70 & 399.24 & 0.05  & 78 \\ \hline
        ~ & 100 & 790.69 & 0.00 & 55 \\ \hline
        ~ & 200 & 1040.18 & 0.00  & 19 \\ \hline
        \rowcolor{Gray}
        50 & ~ & 208.40 & 0.03  &  401 \\ \hline
        ~ & 20 & 0.13 & 0.00 &  90 \\ \hline
        ~ & 30 & 1.91 & 0.01 &  90 \\ \hline
        ~ & 50 & 119.48 & 0.07 &  90 \\ \hline
        ~ & 70 & 456.59 & 0.07 & 72 \\ \hline
        ~ & 100 & 1032.91 & 0.03 & 48 \\ \hline
        ~ & 200 & 728.5 & 11.01 & 11 \\ \hline
        \rowcolor{Gray}
        75 & ~ & 225.11 & 0.14 & 401 \\ \hline
        ~ & 20 & 0.13 & 0.01  & 90 \\ \hline
        ~ & 30 & 1.71 & 0.03 &  90 \\ \hline
        ~ & 50 & 119.28 & 0.27 &  90 \\ \hline
        ~ & 70 & 592.87 & 0.12  &  74 \\ \hline
        ~ & 100 & 1130.55 & 0.48 &  46\\ \hline
        ~ & 200 & 536.80 & 0.6 &  11 \\ \hline
        \rowcolor{Gray}
        100 & ~ & 201.92 & 0.82  & 395 \\ \hline
        ~ & 20 & 0.17 & 0.00  &  90 \\ \hline
        ~ & 30 & 2.13 & 0.14  &  90 \\ \hline
        ~ & 50 & 139.38 & 2.32 &  90  \\ \hline
        ~ & 70 & 665.81 & 0.68  & 73 \\ \hline
        ~ & 100 & 987.60 & 1.35  & 44 \\ \hline
        ~ & 200 & 453.00 & 0.25  & 8 \\ \hline
        \rowcolor{Gray}
        200 & ~ & 208.45 & 43.62  & 329 \\ \hline
        ~ & 20 & 0.14 & 0.34  & 90  \\ \hline
        ~ & 30 & 2.09 & 62.80  & 90  \\ \hline
        ~ & 50 & 184.89 & 144.38 &  90  \\ \hline
        ~ & 70 & 1020.05 & 64.86  &  43 \\ \hline
        ~ & 100 & 552.77 & 0.08  &  13 \\ \hline
        ~ & 200 & 104.00 & 0.00  & 3 \\ \hline
         \rowcolor{Gray}
        ~ & ~ & \textbf{207.30} & \textbf{65.25}  & \textbf{1948} \\ \hline
    \end{tabular}
    \caption{Instances solved by both methods (max-cut)} \label{tMax}
\end{table}

Table \ref{TNonmax} reports the results of the non-concluded instances using the  IP formulation but solved by the algorithm. First of all, the branch and bound algorithm solved these instances in $88.62$ seconds on average. Secondly, the returned solution was coincident in $103$ of the $689$ instances. Regarding the quality of these solutions the $\% GAP$ was $-14.27 \%$. Negative $\% GAP$ means that the IP solution is worse than the one given by the algorithm. Therefore the $\% GAP$ in the solved instances is smaller than in the min-cut case. Finally, although one initial feasible solution was provided to the IP solver, we also have results about what happens when no initial solution was supplied; in thes cases the solver was not even capable of getting feasible solutions for $689$ instances detailed in Column \textit{\# Non-Solution}.

\begin{table}[!ht]
    \centering
    \begin{tabular}{|r|r|r|r|r|r|r|r|}
    \hline
        \textbf{\% p} & \textbf{$\#$ Nodes}  & \textbf{\# Non-concluded} & \textbf{t Algorithm}   & \textbf{$\#$ Coincs.}  & \textbf{\% GAP} & \textbf{\# Non-Solution} \\ \hline
        \rowcolor{Gray}
        25 &    & 124  &  0.34 & 23 & -6.76 & 154\\ \hline
           & 70 & 18 & 0.00  & 11  & 0.00 & 29  \\ \hline
           & 100 & 35 & 0.20 & 7  & -10.59 & 46  \\ \hline
           & 200 & 71 & 1.00  & 5  & -11.27 & 79  \\ \hline   
        \rowcolor{Gray}
        50 &   & 139 &  2.53 & 23   & -8.60 & 173  \\ \hline
           & 70  & 18 & 0.36 & 10  & -5.18 & 32 \\ \hline
           & 100  & 42  & 0.54 & 7  & -14.88 & 55 \\ \hline
           & 200  & 79 & 10.57 & 6  & -3.57 & 86 \\ \hline  
        \rowcolor{Gray}
        75 &    & 176  &  9.92 & 17   & -11.20 & 139 \\ \hline
           & 70  & 30 & 2.43  & 8  & -8.44 & 16  \\ \hline
           & 100  & 61  & 3.29 & 7   & -15.85 & 44 \\ \hline
           & 200  & 85  & 46.17 & 2   & -4.45 & 79 \\ \hline  
        \rowcolor{Gray}
        100 &    & 192 &  101.21  & 21   & -15.43 & 145 \\ \hline
           & 70  & 36 & 30.68 & 8  & -12.83 & 17  \\ \hline
           & 100 & 70 & 22.87 & 11   & -17.60 & 46 \\ \hline
           & 200  & 86 & 8.00  & 2   & -14.41 & 82 \\ \hline  
        \rowcolor{Gray}
        200 &   & 211 &  203.06 & 18   & -21.19  & 142 \\ \hline
           & 70 & 47 & 128.36 & 9    & -17.89 & 17 \\ \hline
           & 100   & 77   & 278.76& 7 & -25.27 & 41 \\ \hline
           & 200 & 87 & 0.00  & 2   & 0.00 & 84  \\ \hline 
        \rowcolor{Gray}
        ~ & ~  & \textbf{906}  & \textbf{88.62} & \textbf{103}  & \textbf{-14.27} & \textbf{689}  \\ 
        \hline
    \end{tabular}
    \caption{Non-solved instances by LP (max-cut)} \label{TNonmax}
\end{table}

The instances not solved by any of the methods are summarized in Table \ref{TNonmax2}. Their $GAP$ was $-71.57 \%$ on average. But for one instance of $70$ nodes and a budget of $200\%$, the IP solution was better than that of the algorithm with a positive GAP of $9.37 \%$.

\begin{table}[!ht]
    \centering
    \begin{tabular}{|r|r|r|r|r|r|r|r|}
    \hline
        \textbf{\% Budget} & \textbf{$\#$ Nodes} & \textbf{\% GAP} & \textbf{\# Non solved} \\ \hline
        75 &  200 & -139.43 & 2\\ \hline
        100 & 200  & -54.59 & 9  \\ \hline 
        100 & 200  & -16.00 & 1  \\ \hline 
        200  & 70 & -39.75 & 5  \\ \hline  
        200  & 70 & 9.37 & 1  \\ \hline
        200  & 100 & -65.49 & 10  \\ \hline  
        200  & 200 & -80.26 & 43  \\ \hline  
        \rowcolor{Gray}
        ~ & ~  & \textbf{-71.57} & \textbf{71}  \\  \hline
    \end{tabular}
    \caption{Non-solved instances by both (max-cut)} \label{TNonmax2}
\end{table}

\section{Conclusions}

In this paper, we have introduced two budget-constrained cut problems: the Budget-Constrained Min-Cut Problem and the Budget-Constrained Max-Cut Problem. \\

We have proposed one IP model for each one of them and also developed a branch and bound enumerative algorithm for solving them. In addition, we have studied some basic properties for these problems. Finally, focusing on the Budget-Constrained Min-Cut Problem, we have also developed a non-exact approach based on a Lagrangian relaxation.\\

In a series of computational experiments, we have also assessed the effectiveness of the different approaches. Clearly, the branch and bound enumerative algorithm outperforms, in terms of computational time, the other methods specially solving the min-cut instances in a reasonable amount of time. The Lagrangian relaxation proposed for the Budget-Constrained Min-Cut Problem has also been very effective although for the considered instances it is  not worth using it since the exact algorithm solves them quickly. We have also analyzed important details about the implementation of the IP formulation in order to get the best possible results.\\

As future work, we consider it interesting to adapt and extend the proposed methods to other types of related cut problems since these algorithms can be seen as a flexible framework for finding constrained cuts. Besides, it is also worth introducing and studying some new problems, such as the k-Min-cut which obtains k minimum cuts of a graph. This last problem could be useful to get exact results from the Lagrangian relaxation that we propose in this paper.\\

\section*{Acknowledgments}
	The first author has been partially supported by the Agencia Estatal de Investigación (AEI) and the European Regional Development Fund (ERDF): PID2020-114594GB--\{C21,C22\}; P18-FR-1422;  and Fundación BBVA: project NetmeetData (Ayudas Fundación BBVA a equipos de investigación científica 2019). The second author has been partially supported by Ministerio de Ciencia e Innovación, project: PID2021-122344NB-I00; and by Generalitat Valenciana, projects: CIGE/2021/161 and CIBEST/2021/35.

\bigskip

\bibliography{budgetCut}

\end{document}